\documentclass[12pt]{amsart}
\usepackage{amsmath}
\usepackage{amsmath}
\usepackage{amsmath}

\newtheorem{theorem}{Theorem}[section]
\newtheorem{lemma}[theorem]{Lemma}
\newtheorem{claim}[theorem]{Claim}
\newtheorem{corollary}[theorem]{Corollary}

\theoremstyle{definition}

\theoremstyle{remark}
\newtheorem{remark}[theorem]{Remark}

\DeclareMathOperator{\Vol}{Vol} \DeclareMathOperator{\dist}{dist}

\DeclareMathOperator{\tr}{tr} \DeclareMathOperator{\inj}{inj}

\begin{document}

\title[Ricci flow on
noncompact manifolds of finite volume]{Non-singular solutions of
normalized Ricci flow on noncompact manifolds of finite volume }

\date{}
\begin{abstract} The main result of this paper shows that, if $g(t)$ is
a complete non-singular solution of the normalized Ricci flow on a
noncompact $4$-manifold $M$ of finite volume, then the Euler
characteristic number $\chi(M)\geq0$. Moreover, $\chi(M)\neq 0$,
there exist a sequence times $t_k\rightarrow\infty$, a double
sequence of points $\{p_{k,l}\}_{l=1}^{N}$ and domains
$\{U_{k,l}\}_{l=1}^{N}$ with $p_{k,l}\in U_{k,l}$ satisfying the
followings:
 \item[(i)] $\dist_{g(t_k)}(p_{k,l_1},p_{k,l_2})\rightarrow\infty$
 as $k\rightarrow\infty$,
 for any fixed $l_1\neq l_2$;
 \item[(ii)] for each $l$, $(U_{k,l},g(t_k),p_{k,l})$ converges in the $C_{loc}^\infty$ sense to a complete
 negative Einstein manifold
 $(M_{\infty,l},g_{\infty,l},p_{\infty,l})$ when
 $k\rightarrow\infty$;
 \item[(iii)]
 $\Vol_{g(t_{k})}(M\backslash\bigcup_{l=1}^{N}U_{k,l})\rightarrow0$ as
 $k\rightarrow\infty$.

\end{abstract}

\author[F. Fang]{Fuquan Fang}
\thanks{The authors were supported by a
NSF Grant of China and the Capital Normal University}
\address{Department of Mathematics, Capital Normal University,
Beijing, P.R.China}
  \email{fuquan\_fang@yahoo.com}
  \author[Z. Zhang]{Zhenlei Zhang}
\address{Department of Mathematics, Capital Normal University,
Beijing, P.R.China  } \email{zhangleigo@yahoo.com.cn}
\author[Y. Zhang]{Yuguang Zhang}
\address{Department of Mathematics, Capital Normal University,
Beijing, P.R.China  }
 \email{zhangyuguang76@yahoo.com}

\maketitle

\section{Introduction}

In his pioneer paper \cite{H99}, Hamilton considered one special
class of Ricci flow solutions on closed three manifolds:
non-singular solutions. Hamilton showed that such solutions provides
an  example of Thurston's geometric decomposition. More precisely,
as the time tends to infinity, the manifolds admit thick-thin
decomposition, where the thick parts converge to hyperbolic spaces,
while the thin parts collapse. In particular,  closed $3$-manifolds
admitting non-singular solutions are geometrizable.

The normalized Ricci flow on a given manifold $M$ is a smooth family
of metrics $g(t),t\in[0,T),$ satisfying the evolution equation
\begin{equation}\label{*}
\frac{\partial}{\partial t}g=-2Ric+\frac{2r}{n}g
\end{equation}
where $Ric$ denotes the Ricci tensor of $g$ and
$r=\frac{\int_MRdv}{\Vol(g)}$ denotes the average scalar curvature
of $g$. The flow requires the scalar curvature to be spatially $L^1$
integrable, so we have to focalize on some special situations.
Following Hamilton \cite{H99}, a solution to equation (\ref{*}) is
called \textit{non-singular} if the solution exists for all time
with uniformly bounded sectional curvature.

In our previous paper \cite{FZZ}, the authors considered
non-singular solutions to the normalized Ricci flow on compact
manifolds and partially generalized Hamilton's convergence results.
Remarkably, the authors found one topological obstruction for the
existence of non-singular solutions on $4$-manifolds, e.g., the
Euler characteristic of the underlying closed $4$-manifold has to be
nonnegative. Moreover, if the Perelman invariant of the $4$-manifold
is negative,  the Hitchin-Thorpe type inequality holds true, i.e.,
$$2\chi (M)\ge \frac 32 |\tau (M)|$$ where $\chi (M)$ (resp. $\tau
(M$) is the Euler characteristic (resp. sinature) of the
$4$-manifold $M$ (cf. \cite{FZZ}). Based on our methods in
\cite{FZZ}, Ishida \cite{Is} recently provided some more examples of
smooth $4$-manifolds which shows that the existence of long time
non-singular solution really depends on the smooth structure of the
underlying manifold.

In this paper we are concerned with non-singular solutions on
complete manifolds with finite volume. Our first result is the
following

\begin{theorem}\label{t11}
Let $g(t)$ be a complete non-singular solution of equation (1) of
finite volume on a noncompact $n$-dimensional manifold $M$, then
either $g(t)$ collapses along a subsequence or $g(t)$ converges
smoothly along a subsequence to a complete Einstein manifold with
negative Einstein constant.
\end{theorem}
Here, "collapse along a subsequence" means that $\max_{x\in
M}\inj(x,g(t_k))\rightarrow0$ for certain sequence of times
$t_k\rightarrow\infty$. The difference from the case of closed
manifolds is the absence of shrinking Ricci solitons and Ricci flat
manifolds in the limit spaces.

In dimension $4$, using the above Theorem 1.1 we obtain the
following

\begin{theorem}\label{t12}
Let $g(t)$ be a complete non-singular solution of equation (1) of
finite volume on a $4$-manifold $M$ (compact or not), then the Euler
characteristic number  $$ 2\chi(M)\geq
 |\frac{1}{16\pi^2}\int_M(|W_{0}^+|^2-|W_{0}^-|^2)dv_{g(0)}|\geq 0,$$ where $W_0^{\pm}$ denotes the
 Weyl tensor  of the initial metric $g(0)$.   Moreover, $\chi(M)=0$ (resp. $\chi (M)=0$ and
 the signature $\text{sig}(M)=0$) if  and only if $g(t)$ collapses along a subsequence
 (resp. $M$ is in addition compact).
\end{theorem}

When $\chi(M)\neq0$, the Ricci flow converges to the negative
Einstein manifolds on the thick part. The volume of the thin part
becomes smaller and smaller, and converges to zero when the time
tends to infinity:
\begin{theorem}\label{t13}
Let $M$ be as in Theorem 1.2.   If $\chi(M)\neq0$, then there exist
a sequence times $t_k\rightarrow\infty$, a double sequence of points
$\{p_{k,l}\}_{l=1}^{N}$ and domains $\{U_{k,l}\}_{l=1}^{N}$ with
$p_{k,l}\in U_{k,l}$ satisfying the followings:
\begin{itemize}
 \item[(i).] $\dist_{g(t_k)}(p_{k,l_1},p_{k,l_2})\rightarrow\infty$
 as $k\rightarrow\infty$,
 for any fixed $l_1\neq l_2$;
 \item[(ii).] for each $l$, $(U_{k,l},g(t_k),p_{k,l})$ converges in the $C_{loc}^\infty$ sense to a complete
 negative Einstein manifold
 $(M_{\infty,l},g_{\infty,l},p_{\infty,l})$ when
 $k\rightarrow\infty$;
 \item[(iii).]
 $\Vol_{g(t_{k})}(M\backslash\bigcup_{l=1}^{N}U_{k,l})\rightarrow0$ as
 $k\rightarrow\infty$.
\end{itemize}
\end{theorem}

Recall that a Riemannian manifold $(M,g)$ is \textit{asymptotic to a
fibred cusp} if $M$ is diffeomorphic to the interior of a compact
manifold $\overline{M}$ whose boundary is a fibration
$F\longrightarrow\partial\overline{M}\longrightarrow B$, and the
metric $g\sim dr^2+\pi^*g_B+e^{-2r}g_F$ at infinity (so the fibres
collapse at infinity). The next theorem concerns the Hitchin-Thorpe
type inequality for non-singular solutions on noncompact four
manifolds asymptotic to a fibered cusp at infinity. We need a
correction term in the Hitchin-Thorpe inequality, namely the
adiabatic limit of the $\eta$-invariants of the infinity. Using the
work of Dai and Wei \cite{DW} we obtain the following theorem:

\begin{theorem}\label{t14}
Let $M$ be as in Theorem 1.1. If $(M,g(0))$ is asymptotic to a
fibred cusp, then the strict Hitchin-Thorpe type inequality holds:
\begin{equation}\label{e11}
2\chi(M)>3|\tau(M)+\frac{1}{2}{\rm a}\lim\eta(\partial\overline{M})|
\end{equation}
where ${\rm a}\lim\eta(\partial\overline{M})$ is the adiabatic limit
of $\eta$ invariant of the boundary.
\end{theorem}

When $\partial\overline{M}$ has special structures, for example
$\partial\overline{M}$ is a disjoint union of circle bundles over
surfaces, we have more precise inequality in the above
Hitchin-Thorpe type inequality (compare \cite{DW}).
\begin{corollary}
Let $(M, g(0))$ be as in Theorem 1.4. If the fibration at infinity
consists of circle bundles over surfaces $S^1\rightarrow
N_i\rightarrow\Sigma_i,1\leq i\leq k,$ then
\begin{equation}
2\chi(M)>3|\tau(M)-\sum_i\frac{1}{3}e_i|
\end{equation}
where $e_i$, $1\le i \le k$, are Euler numbers of the circle
bundles.
\end{corollary}

We remark that the asymptotic assumption in Theorem 1.4 (resp.
Corollary 1.5) may be replaced by assuming the initial metric has
bounded covering geometry and the infinity has a polarized
$F$-structure. Moreover, using the same argument we may extend
Corollary 1.5 to the case where the end of $M$ is asymptotic to a
complex hyperbolic end, i.e,
$$g\sim dr^2+e^{-r}g_{T^2}+e^{-2r}\theta\wedge \theta$$ where $\theta$ is an
invariant $1$-form on the circle fibre, and $\partial \overline{M}$
is a $3$-dimensional nil-manifold.

Comparing with  our previous work in \cite {FZZ}\cite{FZZ2}, it is
natural to ask whether the rigidity theorem could be extended to
noncompact $4$-manifold. More precisely, assume that $M$ is a
complete non-compact $4$-dimensional Riemannian manifold of finite
volume whose end is asymptotic to a complex hyperbolic $4$-manifold,
if $M$ in addition admits a symplectic structure, can one conclude
the Einstein part in Theorem 1.3 is complex hyperbolic under certain
topological constraints? (compare the work of Biquard \cite{Bi}.)

We conclude this introduction by pointing out the main difference
from the compact case dealt in \cite{FZZ}\cite{FZZ2}. To prove the
convergence part of Theorem \ref{t11}, a key lemma we need to verify
in the non-compact case is the vanishing of the integral $\int _M
\triangle R$ (cf. Lemma \ref{l32} below). This follows by using
Shi's derivative estimate for curvatures when the sectional
curvature of the solution is uniformly bounded.  In the proof of
Theorem 1.4 we need to estimate derivatives of the curvature
operator which depends also heavily on Shi's estimation. Because of
this, we really need the sectional curvature bound in the noncompact
case, rather than Ricci bound or even scalar curvature bound in
certain cases as in our previous works. We will get back to this
point in future.

The paper is organized as follows: In Section 2, we recall one
theorem about the maximal principle on noncompact manifolds; in
Section 3, we prove Theorem \ref{t11} and then in Section 4 we prove
Theorem \ref{t12}, Theorem \ref{t13} and Theorem \ref{t14}.

\noindent {\bf Acknowledgement:} The original version of the paper
was written when the second author was visiting UCSD. The second
author would like to thank Professor L. Ni for his invitation. The
second author also would like to thank Prof. L. Ni and B. Chow for
their support and help for living in San Diego.

\section{Preliminaries}

To consider the Ricci flow on noncompact manifolds, we need to use
the maximal principle on noncompact manifolds. For the sake of
reader's convenience let us recall the following general result,
which was proved in \cite{NT}, see also \cite{CLN}.

\begin{theorem}\cite{NT}
Let $(M,g(t)),t\in[0,T],$ be a smooth family of complete evolving
Riemannian manifolds such that $\frac{\partial }{\partial
t}g_{ij}=-2\Upsilon_{ij}$. Denote
$R_\star(t)=\inf_{M}\tr_{g(t)}\Upsilon$ and assume that $R_\star(g)$
is finite and integrable. Assume further that the metrics $g(t)\geq
g^\star$ for a fixed complete metric $g^\star$. Then for any
subsolution to the heat equation $\frac{\partial }{\partial
t}u\leq\triangle_{g(t)}u$, if there is one $\alpha>0$ and $o\in M$
such that
$$\int_0^T\int_M \exp(-\alpha d_\star^2(o,x))u_+^2(x,t)dv_{g(t)}(x)dt<\infty,$$
 where $d_\star$ denotes the distance function of $g^\star$ and
 $u_+=\max(0,u)$, then $u(0)\leq0$ implies $u(t)\leq0$ for all
 time $t\in[0,T]$.
\end{theorem}

One immediate corollary says
\begin{corollary}\cite{CLN}\label{CLN}
Let $(M,g(t)),t\in[0,T],$ be a complete solution to the Ricci flow
with uniformly bounded Ricci curvature. If $u$ is a weak subsolution
of the heat equation on $M\times[0,T]$, such that $u(0)\leq0$ and
$$\int_0^T\int_M\exp(-\alpha d_{g(0)}^2(o,x))u_+^2dv_{g(t)}(x)dt<\infty,$$
for some $\alpha>0$, then $u\leq0$ over $M\times[0,T]$.
\end{corollary}

\section{Non-singular solutions of finite volume}

We will give a proof of Theorem \ref{t11} in this section. Our
argument relies on the following classification theorem of limit
models of Type I Ricci flow, which is due to Naber \cite{Na}. By
Hamilton \cite{H95}, a solution to the Ricci flow
\begin{equation}\label{ricci flow}
\frac{\partial}{\partial t}g(t)=-2Ric(g(t))
\end{equation}is called {\it Type I} if the curvature satisfies
$\sup_{M}|Rm|(t)\leq\frac{C_0}{T-t}$ for some constant $C_0$
independent of $t$.
\begin{theorem}\label{convergence}\cite{Na}
Let $(M,g(t)),t\in[0,T),$ be a Type I solution to the Ricci flow
{\rm(\ref{ricci flow})} such that each metric $g(t)$ is complete,
then for any given point $p\in M$ and times $t_i\rightarrow T$, the
sequence of Ricci flow solutions $(M,(T-t_i)^{-1}g((T-t_i)t+t_i),p)$
will converge along a subsequence to a shrinking Ricci soliton
$(M_\infty,g_\infty(t),p_\infty),t\in(-\infty,1)$.
\end{theorem}

In the following of this section, $M$ stands for a noncompact
manifold and $g(t),t\in[0,\infty),$ is a normalized Ricci flow
solution with finite volume:
\begin{equation}\label{*.5}
\frac{\partial}{\partial t}g=-2Ric+\frac{2r}{n}g,
\end{equation}
where $r=\frac{\int Rdv}{\Vol(g)}$ denotes the average scalar
curvature as usual. Suppose $g(t)$ have uniformly bounded curvatures
$|Rm|(t)\leq C$ for some constant $C$ independent of $t$. Obviously
the flow (6) preserves the volume and so after a scaling we may
assume that $\Vol(g(t))\equiv1$ for all time.

We first establish some lemmas for the proof of the theorem. First
we set $\breve{R}(g)=\inf_M R(g)$ for given metric $g$, the infimum
of the scalar curvature on $M$, and let
$\breve{R}(t)=\breve{R}(g(t))$. The following lemma follows from the
maximal principle.

\begin{lemma}\label{l31} $\breve{R}(t)$ preserves the nonnegative
property and increases whenever it is non-positive.
\end{lemma}
\begin{proof} Consider the evolution
equation of the scalar curvature:
\begin{equation}\label{e6}
\frac{\partial}{\partial t}R=\triangle
R+2|Ric^o|^2+\frac{2}{n}R(R-r),
\end{equation}
where $Ric^o$ denotes the trace free part of the Ricci curvature. We
assume that $r(t)$ is a scalar function defined in advance
satisfying $\breve{R}(t)\leq r(t)$ for all time. Let $f$ be a scalar
function defined by $\frac{d}{dt}f=\frac{2}{n}f(f-r)$ with initial
value $f(0)=\breve{R}(0)-\epsilon$ for fixed small $\epsilon>0$.
Then $f(0)<r(0)$ and $\min(f(0),0)\leq f(t)\leq r(t)$ for all time.
Moreover, $R-f$ satisfies the evolution inequality
$$\frac{\partial}{\partial t}(R-f)\geq\triangle(R-f)+\frac{2}{n}(R+f-r)(R-f),$$
with initial condition $(R-f)(0)>0$ over $M$. Suppose that
$\frac{2}{n}|R+f-r|\leq\bar{C}$, then
$$\frac{\partial}{\partial t}[e^{\bar{C}t}(R-f)]\geq\triangle[e^{\bar{C}t}(R-f)]$$
whenever $R-f\leq0$. Applying Corollary \ref{CLN} to
$e^{\bar{C}t}(R-f)$ we get that $f(t)\leq R(t)$ for all time.

If $R\geq0$ at $t=0$, then $f(0)\geq-\epsilon$ and
$f(t)\geq-\epsilon$ for all time. So $R(t)\geq-\epsilon$ for all
time and then letting $\epsilon$ tend to zero yields the
nonnegativity of $R(t)$. If $\breve{R}(0)\leq0$, then $f(t)\geq
f(0)$ for all time and this gives the monotonicity of $\breve{R}$ by
letting $\epsilon$ tends to zero.
\end{proof}

\begin{lemma}\label{l32}
At each time $t\geq1$, $\int_M\triangle Rdv=0$.
\end{lemma}
\begin{proof}
By Shi's gradient estimate \cite{Sh}, see also \cite{H95}, there is
a constant $C_1<\infty$ such that for any given $t\geq1$, the
estimates $|\triangle R|,|\nabla R|\leq C_1$ spatially holds at time
$t$. Then choose a chopping $\{U_k\}$ with smooth boundaries such
that $U_k\subset U_{k+1}$ and $\bigcup U_k=M$. We can also assume
that $\Vol(\partial U_k,g|_{U_k})\rightarrow0$ as
$k\rightarrow\infty$ because the total volume of $M$ is finite. Then
$$\int_M\triangle Rdv=\lim_{k\rightarrow\infty}\int_{U_k}\triangle Rdv\leq\lim_{k\rightarrow\infty}\int_{\partial U_k}|\nabla R|=0,$$
since $|\nabla R|\leq C_1$ and the volume $\Vol(\partial
U_k)\rightarrow0$.
\end{proof}

To prove the convergence result, we also need the following
\begin{lemma}\label{l33}
If $\breve{R}(t)\leq-c<0$ for all time, then
\begin{eqnarray}
\int_0^\infty(r-\breve{R})dt&<&\infty,\label{e7}\\
\int_0^\infty\int_M|R-r|dvdt<\infty,\label{e8}\\
\int_0^\infty\int_M|Ric^o|^2dvdt&<&\infty.\label{e9}
\end{eqnarray}
\end{lemma}
\begin{proof}
The first estimate follows from the maximal principle. Given
$\epsilon>0$, let $f$ be the function as defined in the proof of
Lemma \ref{l31} which is monotone increasing and satisfies that
$f(t)\leq\breve{R}(t)\leq -c$ for all time. From the evolution
equation $\frac{d}{dt}f=\frac{2}{n}f(f-r)$ we obtain that
$$\int_0^\infty(r(t)-f(t))dt\leq\frac{n}{2c}(-f(0)-c).$$ So
$$\int_0^\infty(r(t)-\breve{R}(t))dt\leq\int_0^\infty(r(t)-f(t))dt<\infty.$$
The second estimate follows directly by
$$\int_0^\infty\int_M|R-r|dvdt\leq\int_0^\infty(R-\breve{R}+r-\breve{R})dvdt=\int_0^\infty2(r-\breve{R})dt.$$

To show the third estimate, we consider the evolution
\begin{eqnarray}
\frac{d}{dt}\int_M Rdv&=&\int_M(\triangle
R+2|Ric^o|^2+\frac{2-n}{n}R(R-r))dv\nonumber\\
&=&\int_M(2|Ric^o|^2+\frac{2-n}{n}R(R-r))dv.\nonumber
\end{eqnarray}
It follows that
\begin{eqnarray}\nonumber
\int_0^\infty\int_M2|Ric^o|^2dvdt\leq\lim_{t\rightarrow\infty}|r(t)-r(0)|+\frac{n-2}{n}\int_0^\infty\int_M
C|R-r|dvdt<\infty,
\end{eqnarray}
which is the desired result.
\end{proof}

The consequence is that the metric tends to be Einstein in the $L^2$
sense:
\begin{lemma}\label{l34}
Suppose as in above lemma, then
\begin{equation}\label{e10}
\lim_{t\rightarrow\infty}(r(t)-\breve{R}(t))=0,
\end{equation}
\begin{equation}\label{e11}
\lim_{t\rightarrow\infty}\int_M|Ric^o|^2dv=0.
\end{equation}
\end{lemma}
\begin{proof}
By above lemma, it suffice to show that
$$\frac{d}{dt}(r-\breve{R})\leq D,\hspace{0.3cm}\mbox{ and }$$
$$\frac{d}{dt}\int_M|Ric^o|^2dv\leq D$$
for some uniform constant $D<\infty$. These facts follow from Shi's
gradient estimate and the non-singular assumption.
\end{proof}

The following is the key lemma for proving Theorem \ref{t11}:
\begin{lemma}\label{l35}
$\liminf_{t\rightarrow\infty}r(t)\leq0$.
\end{lemma}

The proof of this lemma relies on a contradiction argument and we
postpone it to the end of this section.

\begin{lemma}\label{l350} If $\liminf_{t\rightarrow\infty}r(t)=0$,
there is a sequence $t_{k}\longrightarrow \infty$ such that
$$\lim_{k\rightarrow\infty}\int_M2|Ric^o(t_k)|^2dv_{g(t_k)}=0. $$
\end{lemma}
\begin{proof}
Let $t_k\rightarrow\infty$ be a sequence  with
$\lim_{k\rightarrow\infty}r(t_k)=0$. Consider the family of
functions $r(t_k+t)$, $t\in [-t_k, +\infty)$.  By Shi's gradient
estimate \cite{Sh}, for any $l>0$, $$|\frac{d^{l}
r(t_k+t)}{dt^{l}}|\leq C' \sum_{0\leq j \leq l}|\frac{d^{j}
}{dt^{j}}Ric(t_k+t)|\leq C \sum_{0\leq j \leq 2l}|\nabla^{j}
Rm(t_k+t)|\leq \bar{C},
$$ for a constant $\bar{C}>0 $ independent of $t$ and $k$.  By
passing to a subsequence, $r(t_k+t)$ $C^{\infty}$-converges to a
smooth function $r_\infty(t)$ on $\mathbb{R}$, which satisfies
$r_\infty(t)\geq r_\infty(0)=0$.

Indeed, by Eq. (\ref{e10}) in Lemma \ref{l34} and
$\lim_{k\rightarrow\infty}r(t_k)=0$, we deduce that
$\lim_{k\rightarrow\infty}\breve{R}(t_k)=0$. Then by assumption that
$r_\infty(0)=\min_{t\in\mathbb{R}}r_\infty(t)$,
\begin{eqnarray}
0&=&\frac{dr_\infty}{dt}(0)=\lim_{k\rightarrow\infty}\frac{dr}{dt}(t_k)\nonumber\\
&=&\lim_{k\rightarrow\infty}\int_M(2|Ric^o(t_k)|^2+\frac{2-n}{n}R(t_k)(R(t_k)-r(t_k)))dv_{g(t_k)}\nonumber\\
&\geq&\lim_{k\rightarrow\infty}\int_M2|Ric^o(t_k)|^2dv_{g(t_k)}-\lim_{k\rightarrow\infty}\frac{2-n}{n}C\int_M|R(t_k)-r(t_k)|dv_{g(t_k)}\nonumber\\
&\geq&\lim_{k\rightarrow\infty}\int_M2|Ric^o(t_k)|^2dv_{g(t_k)}-\lim_{k\rightarrow\infty}\frac{2-n}{n}C\int_M(R(t_k)+r(t_k)-2\breve{R}(t_k))
dv_{g(t_k)}\nonumber\\
&=&\lim_{k\rightarrow\infty}\int_M2|Ric^o(t_k)|^2dv_{g(t_k)}-\lim_{k\rightarrow\infty}\frac{2-n}{n}2C(r(t_k)-2\breve{R}(t_k))\nonumber\\
&=&\lim_{k\rightarrow\infty}\int_M2|Ric^o(t_k)|^2dv_{g(t_k)}.\nonumber
\end{eqnarray}
\end{proof}

 Now we can give a
\begin{proof}[Proof of Theorem \ref{t11}]
If $\liminf_{t\rightarrow\infty}r(t)=0$, then we claim that for any
sequence of times $t_k\rightarrow\infty$ with
$\lim_{k\rightarrow\infty}r(t_k)=0$, the corresponding metrics
$g(t_k)$ collapse as $k\rightarrow\infty$. Suppose not, then there
exists $\epsilon>0$ such that passing a subsequence, for each $k$,
there is one point $p_k\in M$ with $\inj(p_k,g(t_k))\geq\epsilon>0$.
Then passing one subsequence again, by Hamilton's compactness
theorem \cite{H952}, $(M,g(t_k+t),p_k)$ converge to a limit
"normalized" Ricci flow solution $g_\infty(t)$ on one noncompact
manifold $M_\infty$:
$$\frac{\partial}{\partial
t}g_\infty(t)=-2Ric(g_\infty(t))+\frac{2}{n}r_\infty(t)g_\infty(t),$$
where $r_\infty(t)=\lim_{k\rightarrow\infty}r(t_k+t)$ is nonnegative
satisfying $r_\infty(0)=0$. We will show that $g_\infty(0)$ is Ricci
flat, so $g_\infty(0)$ has infinite total volume (cf. \cite{CGT,
SY}), which contradicts with the fact
$\Vol(M_\infty,g_\infty)\leq1$.

Indeed, by Eq. (\ref{e10}) in Lemma \ref{l34} and
$\lim_{k\rightarrow\infty}r(t_k)=0$, we deduce that
$\lim_{k\rightarrow\infty}\breve{R}(t_k)=0$. Then by Lemma
 \ref{l350}
\begin{eqnarray}
0&=&\lim_{k\rightarrow\infty}\int_M2|Ric^o(t_k)|^2dv_{g(t_k)},\nonumber
\end{eqnarray}
which implies that $g_\infty(0)$ is Einstein.
$\breve{R}(t_k)\rightarrow0$ yields that $g_\infty(0)$ has
nonnegative scalar curvature. Thus $g_\infty(0)$ must be Ricci flat
since the manifold $M_\infty$ is noncompact.

On the other hand, if $\lim\inf_{t\rightarrow\infty}r(t)=-c<0$, then
Eq. (\ref{e10}) implies that
$\lim_{t\rightarrow\infty}\breve{R}(t)=-c<0$. If $g(t)$ do not
collapse along a sequence $t_k\rightarrow\infty$, then $(M,g_{t_k})$
converge subsequently to a limit by Hamilton's compactness theorem
\cite{H952}. By Eq. (\ref{e11}), the limit must be negative
Einstein.
\end{proof}

At last we give a
\begin{proof}[Proof of Lemma \ref{l35}]
Argue by contradiction. By contraries, there is $\delta>0$ such that
$r(t)\geq\delta$ for all time. As showed in \cite{FZZ}, in this
situation, the corresponding unnormalized Ricci flow becomes
singular in finite time and the total volume tends to zero as the
solution approaches the singular time. We can claim more on the
volume decay rate:
\begin{claim}\label{c32}
Let
$\tilde{g}(\tilde{t})=\psi(t)g(t),\tilde{t}\in[0,\widetilde{T}),$ be
the corresponding Ricci flow solution. Then there is $C_2<\infty$
such that
$$C_2^{-1}(\widetilde{T}-\tilde{t})^{n/2}\leq\Vol(\tilde{g}(\tilde{t}))\leq
C_2(\widetilde{T}-\tilde{t})^{n/2}.$$
\end{claim}
\begin{proof}
Comparing the evolution of Ricci flow
\begin{equation}
\frac{\partial}{\partial \tilde{t}}\tilde{g}=-2Ric(\tilde{g}),
\end{equation}
with normalized Ricci flow equation (\ref{*.5}), we obtain the
identities
\begin{eqnarray}
&&\frac{\partial\tilde{t}}{\partial t}=\psi(t) \nonumber\\
&&0=\frac{\partial}{\partial t}(\ln\psi)+\frac{2}{n}r(t).\nonumber
\end{eqnarray}
So the scaling function $\psi(t)=\Vol(\tilde{g}(t))^{2/n}$ is given
by the integration $\psi(t)=\exp(-\int_0^t\frac{2}{n}r(s)ds)$. And
\begin{eqnarray}\nonumber
\tilde{t}=\int_0^{\tilde{t}} d\tilde{t}=\int_0^t\psi(s)ds =\int_0^t
\exp(-\int_0^s\frac{2}{n}r(u)du)ds,
\end{eqnarray}
in particular $\widetilde{T}=\int_0^\infty
\exp(-\int_0^s\frac{2}{n}r(u)du)ds$. Now we can compute
\begin{eqnarray}
(\widetilde{T}-\tilde{t})\Vol(\tilde{g}(\tilde{t}))^{-\frac{2}{n}}&=&\psi(t)^{-1}(\widetilde{T}-\tilde{t})\nonumber\\
&=&\exp(\int_0^t\frac{2}{n}r(u)du)\cdot\int_t^\infty
\exp(-\int_0^s\frac{2}{n}r(u)du)ds\nonumber\\
&=&\int_t^\infty\exp(-\int_t^s\frac{2}{n}r(u)du)ds,\nonumber
\end{eqnarray}
which is bounded from above by $\frac{n}{2\delta}$ and below
$\frac{n}{2C}$. Then letting
$C_2=\max(\frac{n}{2\delta},\frac{2C}{n})$ we get he desired
estimate.
\end{proof}

 Next we claim that
\begin{claim}\label{c33}
The unnormalized Ricci flow
$(M,\tilde{g}(\tilde{t})),\tilde{t}\in[0,\widetilde{T})$, is of Type
I.
\end{claim}
\begin{proof}
By assumption, the average scalar curvature $r(t)$ is comparable
with $\sup|Rm|(t)$ for each time and so it suffice to check that the
quantity
\begin{eqnarray}
\tilde{r}(\tilde{t})(\widetilde{T}-\tilde{t})&=&r(t)\psi(t)^{-1}(\widetilde{T}-\tilde{t})\nonumber
\end{eqnarray}
has a uniform upper bound for all $\tilde{t}<\widetilde{T}$, which
follows directly from the arguments in the proof of Claim \ref{c32}.
\end{proof}

The Claim \ref{c32} shows that
$g(t)=\alpha(\tilde{t})\cdot(\widetilde{T}-\tilde{t})^{-1}\tilde{g}(\tilde{t})$
for certain family of bounded constants
$C_2^{-1}\leq\alpha(\tilde{t})\leq C_2$. Then applying Theorem
\ref{convergence}, $(M,g(t))$ converge along a subsequence to a
noncompact shrinking Ricci soliton, say $(M_\infty,g_\infty)$, with
bounded curvature and finite volume. This is a contradiction because
the volume $\Vol(M_\infty,g_\infty)$ must be infinite. In fact, from
Carrillo and Ni's work \cite{CN}, Perelman's $\mu$ functional
$\mu(g_\infty,1)$ (see \cite{CN} for a definition) is bounded below.
Following Perelman's proof of no local collapsing of finite time
Ricci flow, cf. \cite[\S4]{Pe} or \cite[\S13]{KL}, one can prove
that the volume of any unit metric ball in $M_\infty$ has a uniform
lower bound (notice here the soliton $g_\infty$ has bounded
curvature and so Bishop-Gromov volume comparison theorem works), and
so $\Vol(g_\infty)=\infty$. This finishes the proof of the Lemma
\ref{l35}.
\end{proof}

\vskip 8mm

\section{Non-singular solutions on 4-manifolds with finite volume}
\vskip 4mm

 In this section, we restrict ourselves on
$4$-dimensional case and prove some similar result as in the compact
case (cf. \cite {FZZ}).

\begin{proof}[Proof of Theorem 1.2]
Consider Gauss-Bonnet-Chern formula for a Riemannian manifold of
finite volume and bounded curvature:
\begin{equation}\label{41}
\chi(M)=\frac{1}{8\pi^2}\int_M(\frac{R^2}{24}+|W|^2-\frac{1}{2}|Ric^o|^2)dv.
\end{equation}
First of all, if $g(t)$ collapses along a subsequence, then by
Cheeger-Gromov's $F$-structure theory \cite{CG}, $M$ admits local
non-trivial tori actions and so $\chi(M)=0$. Actually we can say
more about the relation between $\chi(M)$ and the collapsing in our
situation:
\begin{claim}\label{c41}
The following three conditions are equivalent:
\begin{itemize}
 \item[(1)] $\chi(M)=0$; \hspace{0.5cm} {\rm (2)} $\liminf_{t\rightarrow\infty}\breve{R}(t)=0$;
 \item[(3)] $g(t)$ collapses along a subsequence.
 \item[(4)] $\chi(M)=\text{sig}(M)=0$ if $M$ is in addition
 compact.
\end{itemize}
\end{claim}
\begin{proof}[Proof of the Claim]
By above observation and Claim \ref{c32}, it suffice to show that
$(1)$ implies $(2)$ and (3) implies (4). Suppose not, then by Claim
\ref{c32} again, $\lim_{t\rightarrow\infty}\breve{R}(t)=-c$ for some
positive constant $c$. Then by Lemma \ref{l34}, $r(t)\rightarrow-c$
and $\int_M|Ric^o|^2dv\rightarrow0$ as $t\rightarrow\infty$.
Applying the Gauss-Bonnet-Chern formula,
\begin{eqnarray}
\chi(M)&=&\frac{1}{8\pi^2}\int_M(\frac{R^2}{24}+|W|^2-\frac{1}{2}|Ric^o|^2)dv\nonumber\\
&\geq&\frac{1}{192\pi^2}\Vol(g(t))r(t)^2-\frac{1}{16\pi^2}\int_M|Ric^o|^2dv\nonumber\\
&\rightarrow&\frac{c^2}{192\pi^2}\Vol(g(0)),\nonumber
\end{eqnarray}
as $k\rightarrow\infty$, which contradicts with the assumption
$\chi(M)=0$. This finishes the proof of the claim.

To see (3) implies (4) when  $M$ is a compact. We claim that Lemma
\ref{l35} still holds for this case, i.e.
$\liminf_{t\rightarrow\infty}r(t)\leq0$. Otherwise, there is a
sequence $t_{j}\longrightarrow \infty$ such that $g(t_{j}+t) $
collapses, and $r(t_{j}+t)\geq\delta$ for a constant $\delta >0$
independent of $j$ and $t$. However, by the same arguments in the
proof of Lemma \ref{l35}, a subsequence of $g(t_{j}+t) $ converges
to a shrinking soliton $(M_{\infty}, g_{\infty})$, which contradicts
to (3).

By the proof of Lemma \ref{l350} we get a sequence
$t_{k}\longrightarrow \infty$ such that
$$\lim_{k\rightarrow\infty}\int_M2|Ric^o(t_k)|^2dv_{g(t_k)}=0, $$
and, hence (4) follows from below:
 $$2\chi(M)-3|\tau(M)|\geq \lim_{k\rightarrow\infty} \frac{1}{4\pi^2} \int_M (\frac{R^2(t_k)}{24}-\frac{1}{2}|Ric^o(t_k)|^2)dv_{g(t_k)}\geq
 0.$$
\end{proof}

Now we prove that $\chi(M)\geq0$. Indeed, if $\chi(M)\neq0$, then by
above claim and Lemma \ref{l35}, $\breve{R}(t)\rightarrow-c$ for
some $c>0$ and then the previous arguments imply that
$\chi(M)\geq\frac{c^2}{192\pi^2}\Vol(g(0))>0$. The desired result
follows.

Furthermore, we have $$
2\chi(M)-|\frac{1}{16\pi^2}\int_M(|W_{t}^+|^2-|W_{t}^-|^2)dv_{g(t)}|\geq\frac{1}{4\pi^2}\int_M(\frac{R^2}{24}-\frac{1}{2}|Ric^o|^2)dv.$$
By Lemma \ref{l34} and \ref{l350}, there is a sequence
$t_{k}\longrightarrow \infty$ such that $$\lim_{t_{k}\longrightarrow
\infty}\int_M|Ric^o_{t_{k}}|^2dv_{t_{k}}=0,$$ and, hence,
 $$
2\chi(M)\geq \lim_{t_{k}\longrightarrow \infty}
 |\frac{1}{16\pi^2}\int_M(|W_{t_{k}}^+|^2-|W_{t_{k}}^-|^2)dv_{g(t_{k})}|.$$
 Since $(|W_{t}^+|^2-|W_{t}^-|^2)dv_{g(t)}$ is the first  Pontryagin
 form up to a positive  multiplication  by Chern-Weil theory, we
 obtain the desired result from the following claim. By now the
desired result follows from the following claim.\end{proof}

\begin{claim}
For any characteristic polynomial $P$ on $M^4$, we have
\begin{equation}\label{e44}
\int_MP(\Omega_t)=\int_MP(\Omega_0),\hspace{0.3cm}\forall t\geq0.
\end{equation}
\end{claim}
\begin{proof}[Proof of the Claim]The non-singular assumption says that $|Rm(t)|\leq C$ for some
constant $C<\infty$ independent of $t$. Then by Shi's first gradient
estimate \cite{Sh}, there exists another constant $C_4=C_4(n,C)$
such that $|\nabla Rm(t)|\leq C_4(1+t^{-1/2})$ for all $t>0$. Denote
by $\Omega_t$ the curvature operators of metrics $g(t)$.

Given $\rho>0$, denote by $M_\rho$ the set of points whose radial
coordinate at infinity with respect to the metric $g(0)$ are less
that $\rho$. Then by the evolving equation of the volume form, we
have
\begin{equation}\nonumber
|\frac{d}{dt}\ln\Vol_{g(t)}(M\backslash
M_\rho)|=|\Vol_{g(t)}^{-1}(M\backslash M_\rho)\int_{M\backslash
M_\rho}(r-R)dv_{g(t)}|\leq 2C;
\end{equation}
\begin{eqnarray}\nonumber
|\frac{d}{dt}\ln\Vol_{g(t)}(\partial
\overline{M}_\rho)|&=&|\Vol_{g(t)}^{-1}(\partial
\overline{M}_\rho)\int_{\partial
\overline{M}_\rho}(\frac{3}{4}r-R+Ric(\nu,\nu))dv_{g(t)}|\nonumber\\
&\leq&3C,\nonumber
\end{eqnarray}
where $\nu$ is the normal vector fields on $\partial
\overline{M}_\rho$. For fixed time $t_0>0$,
$$\Vol_{g(t_0)}(M\backslash
M_\rho)\leq e^{2Ct_0}\Vol_{g(0)}(M\backslash M_\rho)\rightarrow0$$
and $$\Vol_{g(t_0)}(\partial \overline{M}_\rho)\leq
e^{3Ct_0}\Vol_{g(0)}(\partial \overline{M}_\rho)\rightarrow 0$$ as
$\rho\rightarrow\infty$, since the total volume of $g(0)$ is finite
and $\Vol_{g(0)}(\partial \overline{M}_\rho)\rightarrow0$
exponentially as $\rho\rightarrow\infty$. Combing with Eq.
(\ref{e42}) and (\ref{e43}), and using the curvature bound
$|Rm|(t)\leq C$, it follows that
\begin{eqnarray}
\int_MP(\Omega_{t_0})&=&\lim_{\rho\rightarrow\infty}\int_{M_\rho}P(\Omega_{t_0})\nonumber\\
&=&\lim_{\rho\rightarrow}\int_{M_\rho}P(\Omega_0)
+\lim_{\rho\rightarrow\infty}\int_{M_\rho}(P(\Omega_{t_0})-P(\Omega_0))\nonumber\\
&=&\int_{M}P(\Omega_0)+\lim_{\rho\rightarrow\infty}\int_{M_\rho}2d\int_0^{t_0}P(\dot{\omega}_t,\Omega_t)dt\nonumber\\
&=&\int_MP(\Omega_0)+\lim_{\rho\rightarrow\infty}\int_{\partial
\overline{M}_\rho}\int_0^{t_0}2
P(\dot{\omega}_t,\Omega_t)dt,\label{e45}
\end{eqnarray}
since $deg(P)=2$. In local coordinate $(x^1,\cdots,x^4)$, let
$$\Gamma_{ij}^k(t)=\frac{1}{2}g^{kl}(t)(\frac{\partial g_{il}(t)}{\partial
x^j}+\frac{\partial g_{jl}(t)}{\partial x^i}-\frac{\partial
g_{ij}(t)}{\partial x^l})$$ be the Christoffel symbols of the
Levi-Civita connection at time $t$. Then we can rewrite the
connection one form by $\omega_t(\frac{\partial}{\partial
x^i},\frac{\partial}{\partial
x^j})=\Gamma_{ij}^k(t)\frac{\partial}{\partial x^k}$ and so by the
evolution of the Ricci flow equation (\ref{*}):
$$\dot{\omega}_t(\frac{\partial}{\partial x^i},\frac{\partial}{\partial x^j})=g^{kl}(t)(-\nabla_iR_{jl}-\nabla_jR_{il}+\nabla_lR_{ij})\frac{\partial}{\partial x^{k}},$$
where $\nabla$ denotes the Levi-Civita connection of $g(t)$. Thus
$$|\dot{\omega}_t|\leq3C_4(1+t^{-1/2}).$$ Consequently using the
assumption that $|\Omega_t|=|Rm|\leq C$,
\begin{eqnarray}|\int_{\partial
\overline{M}_{\rho}}\int_0^{t_0}P(\dot{\omega}_t,\Omega_t)dt|
&\leq&C^{''}\int_{\partial \overline{M}_\rho}\int_0^{t_0}|\dot{\omega}_t|\cdot|\Omega_t|dv_{g(t)}dt\nonumber\\
&\leq&3C^{''}C_4C\int_0^{t_0}(1+t^{-1/2})\Vol_{g(t)}(\partial
\overline{M}_\rho)dt\nonumber\\
&\leq&3C^{''}C_4C\Vol_{g(0)}(\partial
\overline{M}_\rho)\cdot\int_0^{t_0}e^{3Ct}(1+t^{-1/2})dt\nonumber\\
&\leq&3C^{''}C_4C(2t_0^{1/2}+t_0)e^{3Ct_0}\cdot\Vol_{g(0)}(\partial
\overline{M}_\rho)\nonumber\\
&\rightarrow&0\nonumber
\end{eqnarray}
as $\rho\rightarrow\infty$, where $C^{''}$ is a constant depending
on the coefficients of the polynomial $P$. Substituting into
(\ref{e45}), the desired result follows.
\end{proof}

\begin{proof}[Proof of Theorem 1.3]
The proof depends on explicit estimate by using Gauss-Bonnet-Chern
formula and Cheeger-Gromov's collapsing theory (with bounded
curvature). We only give a sketch proof here since its proof is
totally the same as in the compact case. By Claim \ref{c41} and
Lemma \ref{l35}, we may assume that $\breve{R}(t)\rightarrow-c<0$ as
$t\rightarrow\infty$.

Choose $t_k\rightarrow\infty$. Fix one small constant
$\varepsilon>0$ such that $M_{k,\epsilon}=:\{x\in
M|\Vol_{g(t_k)}(B_{g(t_k)}(x,1))<\varepsilon\}$ admits an
$F$-structure of positive rank, cf. \cite{CG}. Passing a
subsequence, there is a uniform constant
$N\leq\frac{1}{\varepsilon}$ satisfying that for each $k$, we can
find a maximal set of points $\{p_{k,l}\}_{l=1}^{N_{1}}\subset
M\backslash M_{k,\varepsilon}$ such that as $k\rightarrow\infty$,
$$3\rho_k=\min\{\dist_{g(t_k)}(p_{k,l_1},p_{k,l_2})|l_1\neq
l_2\}\rightarrow\infty,$$
$$\Vol_{g(t_k)}(\bigcup_{l=1}^N\partial B_{g(t_k)}(p_{k,l},\rho_k))\rightarrow0,$$
and the second fundamental forms $\Pi$ of $\partial
B_{g(t_k)}(p_{k,l},\rho_k))$ are uniformly bounded. Then applying
Hamilton's compactness theorem for Ricci flow, cf. \cite{H95}, by
Theorem \ref{t11}, passing a subsequence again we get that for each
$l$, $(B_{g(t_k)}(p_{k,l},\rho_k),g(t_k),p_{k,l})$ converges
smoothly to a complete negative Einstein manifold
$(M_{\infty,l},g_{\infty,l},p_{\infty,l})$.

We next show that $\Vol_{g(t_k)}(M\backslash\bigcup_{l=1}^N
B_{g(t_k)}(p_{k,l},\rho_k))\rightarrow0$ as $k\rightarrow\infty$. By
the choice of the points $\{p_{k,l}\}_{l=1}^N$, we know that
$$M\backslash M_{k,\varepsilon}\subset\{x\in M|\dist(x,\{p_{k,l}\}_{l=1}^N)\leq C_3\}$$
for some constant $C_3$ independent of $k$, when $k$ is large
enough. It concludes that $M\backslash\bigcup_{l=1}^N
B_{g(t_k)}(p_{k,l},\rho_k)\subset M_{k,\varepsilon}$ for large $k$.
So $\chi(M\backslash\bigcup_{l=1}^N B_{g(t_k)}(p_{k,l},\rho_k))=0$
and by Gauss-Bonnet-Chern formula, using the assumption $|R|\leq C$,
\begin{eqnarray}
0&=&\chi(M\backslash\bigcup_{l=1}^N
B_{g(t_k)}(p_{k,l},\rho_k))\nonumber\\
&=&\frac{1}{8\pi^2}\int_{M\backslash\bigcup_{l=1}^N
B_{g(t_k)}(p_{k,l},\rho_k)}(\frac{1}{24}R^2+|W|^2-\frac{1}{2}|Ric^o|^2)dv\nonumber\\
&&+\int_{\bigcup_{l=1}^N\partial
B_{g(t_k)}(p_{k,l},\rho_k)}P(\Pi)dv\nonumber\\
&\geq&\frac{1}{8\pi^2}\int_{M\backslash\bigcup_{l=1}^N
B_{g(t_k)}(p_{k,l},\rho_k)}(\frac{1}{24}\breve{R}^2+\frac{1}{24}(R+\breve{R})(R-\breve{R})-\frac{1}{2}|Ric^o|^2)dv\nonumber\\
&&+\int_{\bigcup_{l=1}^N\partial
B_{g(t_k)}(p_{k,l},\rho_k)}P(\Pi)dv\nonumber\\
&\geq&\frac{c^2}{192\pi^2}\Vol_{g(t_k)}(M\backslash\bigcup_{l=1}^N
B_{g(t_k)}(p_{k,l},\rho_k))\nonumber\\
&&-\frac{1}{8\pi^2}\int_{M\backslash\bigcup_{l=1}^N
B_{g(t_k)}(p_{k,l},\rho_k)}(\frac{C}{12}(R-\breve{R})+\frac{1}{2}|Ric^o|^2)dv\nonumber\\
&&+\int_{\bigcup_{l=1}^N\partial
B_{g(t_k)}(p_{k,l},\rho_k)}P(\Pi)dv\nonumber,
\end{eqnarray}
which implies that $\Vol_{g(t_k)}(M\backslash\bigcup_{l=1}^N
B_{g(t_k)}(p_{k,l},\rho_k))\rightarrow0$ as $k\rightarrow\infty$,
since the last two terms tend to zero by Lemma \ref{l34} and the
assumptions described above. Here $P(\Pi)$ denotes some polynomial
of the second fundamental form $\Pi$. In the last inequality we used
the monotonicity of $\breve{R}$ which implies that $\breve{R}^2\geq
c^2$ for all time. This finishes the proof of the theorem.
\end{proof}

\begin{remark}
By Theorem 1.3 it is easy to see that, if $\chi(M)\neq0$, the number
of Einstein pieces $N\geq1$. Clearly, every piece contributes at
least $1$ to the Euler number $\chi (M)$ and so, $N\le \chi (M)$.
\end{remark}

Before proving Theorem \ref{t14}, let's first recall some groundwork
on the Chern-Weil theory and Chern-Simons correction term, cf.
\cite{Zh}. Let $(N,h)$ be an oriented Riemannian $2n$-manifold. By
Chern-Weil theory, any $SO(2n)$ invariant polynomial of degree $n$,
say $P$, defines a characteristic form $P(\Omega)$, where
$\Omega\in\Lambda^2N\otimes\Lambda^2N$ denotes the curvature
operator. If we have a smooth family of metrics $h_t,t\in[t_1,t_2],$
then the Chern-Simons form $Q_P$, associated to $P$, is defined by
the equation:
\begin{equation}\label{e42}
Q_P(h_{t_2},h_{t_1})=n\int_{t_1}^{t_2}P(\dot{\omega}_t,\Omega_t,\cdots,\Omega_t)dt
\end{equation}
which determines the nice correction term
\begin{equation}\label{e43}
P(\Omega_{t_2})-P(\Omega_{t_1})=dQ_P(h_{t_2},h_{t_1}),
\end{equation}
where $\omega_t$ and $\Omega_t$ denote the connection one form and
the curvature form of the metric $h_t$ respectively. In our
consideration, $P$ will be the Pfaffian $Pf$ or $L$-polynomial
characteristic form.

Using Claim 4.2, and  the Atiyah-Patodi-Singer index formula
\cite{APS} on manifolds with boundary, Dai and Wei proved in
\cite{DW} the following theorem for manifolds with fibred cuspidal
infinity:
\begin{theorem}\cite{DW}\label{t43}
Let $(N,h)$ be a complete Riemannian 4-manifold which is asymptotic
to a fibred cusp metric at infinity, then the Euler number $\chi(N)$
and signature $\tau(N)$ are given by
\begin{equation}
\chi(N)=\int_NPf(\frac{\Omega}{2\pi});
\end{equation}
\begin{equation}
\tau(N)=\int_NL(\frac{\Omega}{2\pi})-\frac{1}{2}{\rm
a}\lim\eta(\partial N),
\end{equation}
where ${\rm a}\lim\eta(\partial N)$ denotes the adiabatic limit of
$\eta(\partial N)$.
\end{theorem}

Here the adiabatic limit ${\rm a}\lim\eta(\partial N)$ is a
topological invariant of the 3-manifold $\partial N$. Now we are
ready to give a proof of Theorem \ref{t14}:
\begin{proof}[Proof of Theorem \ref{t14}]
Let $g(t)$ be a non-singular solution on noncompact $4$-manifold $M$
such that $g(0)$ has asymptotical fibred cusps at infinity.
Topologically $M$ is the interior of a manifold $\overline{M}$ whose
boundary admits a fibration structure
$$F\longrightarrow\partial\overline{M}\stackrel{\pi}{\longrightarrow}B$$
for closed manifolds $B,F$ and geometrically the metric $g(0)$ has
the form
$$g(0)\sim dr^2+\pi^*g_B+e^{-2r}g_F$$
at infinity, where $g_B$ is a metric on $B$ and $g_F=g_F(b)$, $b\in
B$, is a family of metrics on $F$.

Using Dai and Wei's Theorem \ref{t43}, we obtain that at any time
$t$,
\begin{equation}\label{e46}
\chi(M)=\int_MPf(\Omega_0)=\int_MPf(\Omega_t);
\end{equation}
\begin{equation}\label{e47}
\tau(M)=\int_ML(\Omega_0)-\frac{1}{2}{\rm a}\lim\eta(\partial
\overline{M})=\int_ML(\Omega_t)-\frac{1}{2}{\rm a}\lim\eta(\partial
\overline{M}).
\end{equation}
More precisely, at any time $t$, we have the Gauss-Bonnet-Chern
formula
\begin{equation}\label{e49}
\chi(M)=\frac{1}{8\pi^2}\int_M(\frac{R^2}{24}+\frac{1}{4}|W|^2-\frac{1}{2}|Ric^o|^2)dv_{g(t)}
\end{equation}
and the generalized Hirzebruch signature formula
\begin{equation}\label{e49}
\tau(M)=\frac{1}{48\pi^2}\int_M(|W^+|^2-|W^-|^2)dv_{g(t)}-\frac{1}{2}{\rm
a}\lim\eta(\partial \overline{M}).
\end{equation}
It follows that at any time $t$,
\begin{equation}\label{410}
2\chi(M)-2|\tau(M)+\frac{1}{2}{\rm a}\lim\eta(\partial
\overline{M})|\geq\frac{1}{4\pi^2}\int_M(\frac{R^2}{24}-\frac{1}{2}|Ric^o|^2)dv.
\end{equation}
Then combining Lemma \ref{l35}, Claim \ref{c41} and Lemma \ref{l34}
derives the desired strict Hitchin-Thorpe type inequality
(\ref{e11}) by letting $t\rightarrow\infty$.
\end{proof}

\end{document}